\documentclass{amsart}
\usepackage{a4wide}
\usepackage[backref=page,bookmarksopen=true]{hyperref}
\usepackage[latin1]{inputenc}  % accents 8 bits dans le source
\usepackage[T1]{fontenc}       % accents dans le DVI
\usepackage[english]{babel}
\usepackage{amssymb}
\usepackage[alphabetic]{amsrefs}
\usepackage{latexsym}
\usepackage{mathrsfs}
\usepackage{graphics, graphicx}
\usepackage[all, cmtip]{xy}
\usepackage{tabularx}
\usepackage{multirow}
\usepackage{longtable}

\theoremstyle{plain}
\newtheorem{prop}{Proposition}[section]

\newtheorem{cor}[prop]{Corollary}
\newtheorem{lem}[prop]{Lemma}
\newtheorem{thmintro}{Theorem}
\renewcommand{\thethmintro}{\Alph{thmintro}}
\newtheorem{corintro}[thmintro]{Corollary}

\theoremstyle{definition}

\newtheorem*{rem}{Remark}
\theoremstyle{remark}

\newcommand{\Saff}{\mathcal{S}_{\mathrm{aff}}}

\newcommand{\Z}{\mathbb{Z}}

\newcommand{\la}{\langle}
\newcommand{\ra}{\rangle}
\newcommand{\inv}{^{-1}}
\newcommand{\MF}{\mathscr{M}(F)}
\newcommand{\Meucl}{\mathscr{M}_{\mathrm{Eucl}}(F)}

\newcommand{\Iaff}{I_{\mathrm{aff}}}

\newcommand{\minus}{\backslash}
\newcommand{\CAT}[1]{\mathrm{CAT}(#1)}%
\makeatletter

\@addtoreset{equation}{section} \makeatother
\DeclareMathOperator{\Pc}{Pc}
\DeclareMathOperator{\Isom}{Isom}
\DeclareMathOperator{\Ch}{Ch}

\DeclareMathOperator{\diam}{diam}
\newcommand{\norma}{\mathscr{N}}
\newcommand{\centra}{\mathscr{Z}}

\long\def\symbolfootnote[#1]#2{\begingroup%
\def\thefootnote{\fnsymbol{footnote}}\footnote[#1]{#2}\endgroup}
\begin{document}
\title[Relatively hyperbolic  Coxeter groups]
{Buildings with isolated subspaces and\\ relatively hyperbolic Coxeter groups}
\author[P.-E. Caprace]{Pierre-Emmanuel Caprace$^*$}
\thanks{$^*$F.N.R.S. research associate. This work was partially supported by IPDE}
\address{Universit\'e catholique de Louvain, Chemin du Cyclotron 2, 1348 Louvain-la-Neuve, Belgium.}
\email{pe.caprace@uclouvain.be}

\date{First draft, March 2007. Revised version, June 2008. Erratum appended in December 2013}

\begin{abstract}
Let $(W, S)$ be a Coxeter system. We give necessary and sufficient conditions on the Coxeter diagram of $(W, S)$
for $W$ to be relatively hyperbolic with respect to a collection of finitely generated subgroups. The peripheral
subgroups are necessarily parabolic subgroups (in the sense of Coxeter group theory). As an application, we
present a criterion for the maximal flats of the Davis complex of $(W,S)$ to be isolated. If this is the case,
then the maximal affine sub-buildings of any building of type $(W,S)$ are isolated.
\end{abstract}

%\newcommand{\subjclass}[1]{\symbolfootnote[0]{\noindent AMS classification numbers (2000):~#1.}}
%\newcommand{\keywords}[1]{\symbolfootnote[0]{\emph{Keywords}:~#1}}

%\subjclass{20F67; 20E42, 20F55, 20F69}%% AMS classification numbers
\keywords{Coxeter group, building, isolated flat, relative hyperbolicity.}

\maketitle

\section{Introduction}

Let $X$ be a complete $\CAT0$ space. A $k$--\textbf{flat} in $X$ is a subset which is isometric to the
$k$--dimensional Euclidean space. Since we will mainly be interested in isolated flats, it is convenient to
define a \textbf{flat} as a $k$--flat for some $k \geq 2$. In particular, geodesic lines are not considered to
be flats. Let $\mathscr{F}$ be a collection of closed convex subsets of $X$. We say that the elements of
$\mathscr{F}$ are \textbf{isolated} in $X$ if the following conditions hold:
\begin{description}
\item[(A)] There is a constant $D<\infty$ such that each flat $F$ of $X$ lies in a $D$--tubular neighbourhood of
some $C\in \mathscr{F}$.

\item[(B)] For each positive $r<\infty$ there is a constant $\rho=\rho(r)<\infty$ so that for any two distinct
elements $C, C' \in \mathscr{F}$ we have
$$
   \diam \bigl( \mathcal{N}_r(C) \cap \mathcal{N}_r(C') \bigr) < \rho,
$$
where $\mathcal{N}_r(C)$ denotes the $r$--tubular neighbourhood of $C$.
\end{description}
We say that $X$ has \textbf{isolated flats} if $\mathscr{F}$ consists of flats.

Let now $(W,S)$ be a Coxeter system with $S$ finite. Given a subset $J \subset S$, we set $W_J = \la J \ra$; the
group $W_J$ as well as any of its $W$-conjugate, is called a \textbf{parabolic subgroup} of $W$. We also set
$J^\perp = (S \minus J) \cap  \centra_W(W_J)$. The set  $J$ is called \textbf{spherical} (resp.
\textbf{irreducible affine}, \textbf{affine}, \textbf{Euclidean}) if $W_J$ is finite (resp. an irreducible
affine Coxeter group, a direct product of irreducible affine Coxeter groups, a direct product of finite and
affine Coxeter groups). We say that $J$ is \textbf{minimal hyperbolic} if it is non-spherical and non-affine but
every proper subset is spherical or irreducible affine. Let $X(W, S)$ be the Davis complex of $(W, S)$. Thus
$X(W, S)$ is a proper $\CAT0$ space \cite{Da98} and its isometry group contains $W$ as a cocompact lattice.

According to a theorem of G.~Moussong \cite{Moussong}, the group $W$ is Gromov hyperbolic if and only if $S$ has
no irreducible affine subset of cardinality $\geq 3$ and if for each non-spherical $J \subset S$, the set
$J^\perp$ is spherical. Our main result gives necessary and sufficient conditions for $W$ to be relatively
hyperbolic with respect to a collection of parabolic subgroups:

\begin{thmintro}\label{thm:RH}
Let $(W, S)$ be a Coxeter system with $S$ finite, let $\mathscr{P}$ be a collection of parabolic subgroups of
$W$ and let $\mathscr{T}$ be the set of types of elements of $\mathscr{P}$. Then the following conditions are
equivalent:
\begin{enumerate}
\item $\mathscr{T}$ satisfies the following conditions:
\begin{description}
\item[(RH1)] For each irreducible affine subset $J \subset S$, there exists $K \in \mathscr{T}$ such that  $J
\subset K$. Similarly, for each pair of irreducible non-spherical subsets $J_1, J_2 \subset S$ with $[J_1, J_2]
= 1$, there exists $K \in \mathscr{T}$ such that $J_1 \cup J_2 \subset K$.

\item[(RH2)] For all $K_1, K_2 \in \mathscr{T}$ with $K_1 \neq K_2$, the intersection $K_1 \cap K_2$ is
spherical.
\end{description}

\item $W$ is relatively hyperbolic with respect to $\mathscr{P}$.

\item In the Davis complex $X(W, S)$, the residues whose type belongs to $\mathscr{T}$ are isolated.

\item In any building of type $(W, S)$, the residues whose type belongs to $\mathscr{T}$ are isolated.
\end{enumerate}
\end{thmintro}

Basic definitions and properties of relatively hyperbolic groups may be consulted in the standard references \cite{Bowditch_relhyp} or \cite{Farb}. 

\begin{rem}
Throughout this paper, the term `\emph{parabolic subgroup}' will be used only in the sense which was defined above; this agrees with the standard conventions in the theory of Coxeter groups. Given a group $G$ which is relatively hyperbolic with respect to subgroups $H_1, \dots, H_n$, it is customary to call any conjugate of some $H_i$ a \emph{parabolic subgroup} of $G$. In order to avoid any confusion in the present paper, we shall instead call these the \textbf{peripheral subgroups} of $G$. Thus the term `\emph{parabolic}' will be exclusively used in its Coxeter group acceptation.
\end{rem}

Notice that we do not assume the buildings to be locally compact in (iv). Conditions (RH1) and (RH2) can be
checked concretely on the Coxeter diagram of $(W, S)$. Combining Theorem~\ref{thm:RH} with the following, one
obtains in particular a complete characterization of those Coxeter groups which are relatively hyperbolic with
respect to any family of finitely generated subgroups:

\begin{thmintro}\label{thm:periph}
Let $(W,S)$ be a Coxeter system with $S$ finite. If $W$ is relatively hyperbolic with respect to finitely
generated subgroups $H_1, \dots, H_m$, then each $H_i$ is a parabolic subgroup of $W$.
\end{thmintro}

It should be noted that there exist non-affine Coxeter groups which are not relatively hyperbolic with respect
to any family of parabolic subgroups. Consider for example the Coxeter group $W$ with Coxeter generating set $S
= \{s_1, \dots, s_n\}$ defined by the following relations: $[s_i, s_j] = 1$ for $|i-j| \geq 2$ and $o(s_i s_j )
= 4$ for $|i-j|=1$. It is easily verified, using Theorem~\ref{thm:RH}, that for $n > 7$, the group $W$ is not
relatively hyperbolic with respect to any collection of proper parabolic subgroups. For $n = 7$, one checks that
the set
$$\mathscr{T} = \Big\{ \{1, 2, 3, 5, 6, 7\}, \{2, 3, 4\}, \{3, 4, 5\}, \{4, 5, 6\} \Big\}$$
satisfies (RH1) and (RH2).

It happens however quite often that a Coxeter group is relatively hyperbolic with respect to a maximal proper
parabolic subgroup:

\begin{corintro}\label{cor:MaxParab}
Suppose that there exists an element $s_0 \in S$ such that $\{s_0\}^\perp$ is spherical. Then $W$ is relatively
hyperbolic with respect to the parabolic subgroups whose type belongs to the set
$$
\mathscr{T} = \Big\{S\minus \{s_0\}\Big\} \cup \Big\{ J \subset S \; | \; \; J \text{ is affine and contains}\
s_0\Big\}.
$$
\end{corintro}

Of particular interest is the special case when the peripheral subgroups of a relatively hyperbolic group are
virtually nilpotent (or more generally amenable). As for any discrete group acting properly and cocompactly on a
$\CAT0$ space, any amenable subgroup of a Coxeter group is virtually abelian. Theorem~\ref{thm:RH} yields the
following characterization:

\begin{corintro}\label{thm:A}
  The following assertions are equivalent:
  \begin{itemize}
    \item[(i)] For all non-spherical $J_1, J_2 \subset S$ such that $[J_1, J_2]=1$, the group $\la J_1 \cup J_2 \ra$
    is virtually abelian.

    \item[(ii)] For each minimal hyperbolic $J \subset S$, the set $J^\perp$ is spherical.

    \item[(iii)] The Davis complex $X(W, S)$ has isolated flats.

    \item[(iv)] The group $W$ is relatively hyperbolic with respect to a collection of virtually abelian
    subgroups of rank at least~$2$.

    \item[(v)] In any building of type $(W, S)$, the maximal residues of (non-spherical)
    Euclidean type (i.e. the maximal non-spherical Euclidean sub-buildings) are isolated.
  \end{itemize}
\end{corintro}

The list of all minimal hyperbolic Coxeter systems is known (see  \cite[Ch.~V, pp.133--134]{Bo68} or \cite[II.6.9]{Humphreys}); in fact, it is not difficult to see
that a minimal hyperbolic subset of $S$ has cardinality at most~$10$. Many Coxeter systems $(W,S)$ such that $W$
is not Gromov hyperbolic do satisfy condition (ii). In particular $S$ may contain affine subsets.

Theorem~\ref{thm:RH} is deduced from the detailed study of flats in buildings and Coxeter groups which is made
in \cite{CH06}. The equivalence between (ii) and (iii) is a consequence of
\cite[Appendix]{HruskaKleiner_IsolatedFlats}. In fact, the theorem above allows one to apply varied algebraic
and geometric consequences of the isolation of subspaces established in [loc. cit.] and \cite{DrutuSapir}. We
collect a few of them in the special case of virtually abelian peripheral subgroups:

\begin{corintro}\label{cor:B}
Assume that $(W, S)$ satisfies the equivalent conditions of Corollary~\ref{thm:A}. Let $Y$ be a building of type
$(W, S)$, $\mathscr{F}$ be the collection of maximal residues of non-spherical Euclidean type and $\Gamma <
\Isom(Y)$ be a subgroup acting properly discontinuously and cocompactly. Then:
\begin{itemize}
  \item[(i)] $\Gamma$ is relatively hyperbolic with respect to the family of stabilizers of elements of
  $\mathscr{F}$; each of these stabilizers is a cocompact lattice of a Euclidean building.

  \item[(ii)] $W$ and $\Gamma$ are biautomatic.

  \item[(iii)] Every connected component of $\partial_T X(W, S)$ (resp. $Y$) is either an isolated point or a Euclidean
  sphere (resp. a spherical building).

  \item[(iv)] Every asymptotic cone of $X(W, S)$ (resp. $Y$) is tree-graded with respect to a family of
  closed convex subsets which are flats (resp. Euclidean buildings);
  furthermore any quasi-isometry of $X(W, S)$ (resp. $Y$) permutes these pieces.
\end{itemize}
\end{corintro}

We refer to \cite{DrutuSapir} for more information on asymptotic cones and tree-graded spaces. It is known that
all Coxeter groups are automatic \cite{BH94}, but the problem of determining which Coxeter groups are
biautomatic is still incompletely solved: it follows from \cite{CM05} that $W$ is biautomatic whenever $S$ has
no irreducible affine subset of cardinality $\geq 3$. Corollary~\ref{cor:B} shows biautomaticity in many other
cases.

Let us finally mention that the construction of cocompact lattices in $\Isom(Y)$ is a delicate problem, unless
the Coxeter system $(W, S)$ is \textbf{right-angled} (i.e. $o(st) \in \{1, 2, \infty\}$ for all $s, t \in S$).
This question seems especially interesting when $(W, S)$ is \textbf{$2$--spherical}, namely $o(st) < \infty$ for
all $s, t \in S$. Besides the classical case of Euclidean buildings, some known constructions provide examples
of lattices when $W$ is $2$--spherical and Gromov hyperbolic \cite{KangVdovina06}. However, I don't know any
example of a cocompact lattice in $\Isom(Y)$ in the case when $W$ is a $2$--spherical Coxeter group which is
neither Euclidean nor Gromov--hyperbolic. The situation is completely different when $(W, S)$ is right-angled.
In that case indeed, graphs of groups provide a large family of examples of cocompact lattices  to which
Corollary~\ref{cor:B} may be applied.

In order to state this properly, let $A$ be a finite simple graph with vertex set $I$ and for each $i \in I$,
let $P_i$ be a group. Let $\Gamma = \Gamma(A, (P_i)_{i \in I})$ be the group which is the quotient of the free
product of the $(P_i)_{i \in I}$ by the normal subgroup generated by all commutators of the form $[g_i, g_j]$
with $g_i \in P_i$, $g_j \in P_j$ and $\{i, j\}$ spanning an edge of $A$. Let also $(W, \{s_i\}_{i \in I})$ be
the Coxeter system such that $o(s_i s_j)=2$ (resp. $o(s_i s_j)=\infty$) for each edge (resp. non-edge) $\{i,
j\}$ of $A$. Then $\Gamma$ acts simply transitively on the chambers of a building $Y(A, (P_i)_{i \in I})$ of
type $(W, \{s_i\}_{i \in I})$ by \cite[Theorem~5.1 and Corollary~11.7]{Da98}. If each $P_i$ is finite, then this
building is locally compact and, hence, $\Gamma$ is a cocompact lattice in its automorphism group. For example,
if the graph $A$ is a $n$-cycle with $n \geq 5$, then this building is a Bourdon building and $\Gamma$ is a
Bourdon lattice (these are the Fuchsian buildings and their lattices defined and studied by M.~Bourdon in \cite{Bourdon}). Moreover, if each $P_i$ is infinite cyclic, then $\Gamma$ is a right-angled Artin group.

For the Coxeter system $(W, \{s_i\}_{i \in I})$, condition (ii) of Corollary~\ref{thm:A} may be expressed as
follows: \emph{for each $3$-subset $J \subset I$ which is not a triangle, the subgraph induced on $J^\perp = \{i
\in I \; | \; $\{i, j\}$ \text{ is an edge for each } j \in J\}$ is a complete graph}. Let now $\Iaff$ be the
set of all subsets of $J$ of the form $\{i_1, j_1, \dots, i_n, j_n\}$ where $\{i_k, j_k\}$ is a non-edge for
each $k$ but all other pairs of elements are edges. Note that for $J =\{i_1, j_1, \dots, i_n, j_n\}$, the
subgroup $\Gamma_J$ of $\Gamma$ generated by all $P_i$'s with $i \in J$ has the following structure:
$$\Gamma_J \simeq (P_{i_1} * P_{j_1}) \times \dots \times (P_{i_n} * P_{j_n}).$$
Then Corollary~\ref{cor:B} implies that, under the assumption (ii), the group $\Gamma$ is relatively hyperbolic
with respect to the family of all conjugates of subgroups of the form $\Gamma_J$ with $J \in \Saff$.

\section{On parabolic subgroups of Coxeter groups}

Recall that a subgroup of $W$ of the form $W_J$ for some $J \subset S$ is called a \textbf{standard parabolic
subgroup}. Any of its conjugates is called a \textbf{parabolic subgroup} of $W$. A basic fact on Coxeter groups
is that any intersection of parabolic subgroups is itself a parabolic subgroup. This allows to define the
\textbf{parabolic closure} $\Pc(R)$ of a subset $R \subset W$: it is the smallest parabolic subgroup of $W$
containing $R$.

\begin{lem}\label{lem:Pc}
  Let $G$ be a reflection subgroup of $W$ (i.e. a subgroup of $W$ generated by reflections). Then there is a set
  of reflections $R \subset G$ such that $(G, R)$ is a Coxeter system. Furthermore, if $(G, R)$ is irreducible
  (resp. spherical, affine of rank $\geq 3$), then so is $\Pc(G)$.
\end{lem}
\begin{proof}
  For the first assertion, see \cite{Deo89}. Any two reflections in $R$ which do not commute lie in the same
  irreducible component of $\Pc(R)$. Therefore, if $(G, R)$ is irreducible, then all elements of $R$ are in the
  same irreducible component of $\Pc(R)$. Since $G = \la R \ra$ and $\Pc(R)$ is the minimal parabolic subgroup
  containing $G$, we deduce that $\Pc(R)$ is irreducible. If $G$ is finite, then it is contained in a finite
  parabolic subgroup (see \cite{Bo68}), %find Ref Exercise???
  hence $\Pc(G)$ is spherical. Finally, if $(G, R)$ is affine of rank $\geq
  3$, then so is $\Pc(G)$ by \cite[Proposition~16]{Ca06}.
\end{proof}

\begin{lem}\label{lem:Normalizer}
Let $P \subset W$ be an infinite irreducible parabolic subgroup. Then the normalizer of $P$ in $W$ splits as a
direct product: $\norma_W(P) = P \times \centra_W(P)$ and $\centra_W(P)$ is also a parabolic subgroup of $W$.
\end{lem}
\begin{proof}
See \cite[Proposition~5.5]{Deodhar82}.
\end{proof}

\begin{lem}\label{lem:DirectProd}
  Let $G_1, G_2$ be reflection subgroups of $W$ which are irreducible, i.e. $(G_i, R_i)$ is
  irreducible for $R_i \subset G_i$ as in Lemma~\ref{lem:Pc}, and assume that $G_1$ is infinite and that
  $[G_1, G_2] = \{1\}$. Then either
  $$\Pc(\la G_1 \cup G_2 \ra) \simeq \Pc(G_1) \times \Pc(G_2)$$
  or $\Pc(G_1) = \Pc(G_2)$ is an irreducible affine Coxeter group of rank $\geq 3$.
\end{lem}
\begin{proof}
By Lemma~\ref{lem:Pc}, the parabolic closure $\Pc(G_1)$ is infinite and irreducible. Given a reflection $r \in
G_2$, then $r$ centralizes $G_1$ by hypothesis, hence it normalizes $\Pc(G_1)$. Thus either $r
\in \Pc(G_1)$ or $r$ centralizes $\Pc(G_1)$ by Lemma~\ref{lem:Normalizer}. Therefore, either $G_2 \subset
\Pc(G_1)$ or $G_2$ centralizes $\Pc(G_1)$ and $G_2 \cap \Pc(G_1) =\{1\}$.

If $G_2$ centralizes $\Pc(G_1)$ and $G_2 \cap \Pc(G_1) =\{1\}$, then $\Pc(G_2)$ centralizes $\Pc(G_1)$ and
$\Pc(G_1) \cap \Pc(G_2)=\{1\}$ by Lemma~\ref{lem:Normalizer}. Hence we are done in this case.

Assume now that $G_2 \subset \Pc(G_1)$. Then, since $G_1$ normalizes $\Pc(G_2)$, we deduce from
Lemma~\ref{lem:Normalizer} that $\Pc(G_1)=\Pc(G_2)$. It is well known and not difficult to see that $G_1$
contains an element of infinite order $w_1$ such that $\Pc(w_1) = \Pc(G_1)$ (take for example $w_1$ to be the
Coxeter element in the Coxeter system $(G_1, R_1)$ provided by Lemma~\ref{lem:Pc}). Similarly, let $w_2 \in G_2$ be such that $\Pc(w_2) = \Pc(G_2) = \Pc(G_1)$. Thus $w_1$ and $w_2$ are
mutually centralizing. Moreover, we have $\la w_1 \ra \cap \la w_2 \ra < G_1 \cap G_2  = \{1\}$: indeed, any infinite irreducible Coxeter group  is center-free by
Lemma~\ref{lem:Normalizer}. Thus $\la w_1, w_2 \ra \simeq \Z \times \Z$. By \cite[Corollary~6.3.10]{Kra94}, this
implies that $\Pc(G_1) = \Pc(G_2)$ is affine, and clearly of rank $\geq 3$ since it contains $\Z \times \Z$.
\end{proof}

Let $X = X(W, S)$ denote the Davis complex.

\begin{lem}\label{lem:debile}
Let $r, r', s, t$ be reflections. Assume that the wall $X^t$ separates $X^r$ from $X^{r'}$ and that $s$ commutes
with both $r$ and $r'$. Then either $s$ also commutes with $t$ or $s$ belongs to the parabolic closure of $\la
r, r'\ra$.
\end{lem}
\begin{proof}
Let $H < W$ be the infinite dihedral subgroup generated by $r$ and $r'$. By assumption $s$ centralizes $H$,
whence $s$ normalizes the parabolic $\Pc(H)$. By Lemma~\ref{lem:Pc}, the parabolic subgroup $\Pc(H)$ is
irreducible and non-spherical. Since $s$ is a reflection, we deduce from Lemma~\ref{lem:Normalizer} that either
$s $ belongs to $\Pc(H)$ or $s$ centralizes $\Pc(H)$. This finishes the proof because, by \cite[Lemma~17]{Ca06},
the reflection $t$ belongs to $\Pc(H)$.
\end{proof}

\section{On Euclidean flats in the $\CAT0$ realization of Tits buildings}

Let now $F$ be a flat in $X = X(W, S)$; we remind the reader that $\dim(F) \geq 2$ according to the convention adopted in this paper. We use the notation and terminology of \cite{CH06}. In particular, we
denote by $\mathscr{M}(F)$ the set of all walls which separate points of $F$. Furthermore, for any set of walls
$M$, we denote by $W(M)$ the subgroup of $W$ generated by all reflections through walls in $M$. For any $m \in
\MF$, the set $m \cap F$ is a Euclidean hyperplane of $F$ \cite[Lemma~4.1]{CH06}. Two elements $m, m'$ of $\MF$
are called \textbf{$F$--parallel} if the hyperplanes $m \cap F$ and $m' \cap F$ are parallel in $F$. The
following result collects some key facts on flats in $X$ established in \cite{CH06}:

\begin{prop}\label{prop}
  Let $F$ be a flat in $X = X(W, S)$. Then the $F$--parallelism in $\MF$ induces a partition
  $$\MF = M_0 \cup M_1 \cup \dots \cup M_k$$%
  such that
  %
  %$$W(\MF) \simeq W(M_0) \times \dots \times W(M_k),$$%
$$\Pc(W(\MF)) \simeq \Pc(W(M_0)) \times \dots \times \Pc(W(M_k)),$$
  where each $W(M_i)$ is a direct product of infinite irreducible Coxeter groups and
  $W(M_0)$ is a direct product of irreducible
  affine Coxeter groups. Moreover, the set $M_i$ is non-empty for each $i \geq 1$ and
  if $M_0 = \varnothing$, then $k \geq \dim F$.
%  Furthermore, we have
%  %
%  $$\Pc(W(\MF)) \simeq \Pc(W(M_0)) \times \dots \times \Pc(W(M_k)).$$
%  %
\end{prop}
\begin{proof}
Let $\MF = M_1 \cup \dots \cup M_l$ be the partition of $\MF$ into $F$--parallelism classes. Since the dimension
of $F$ is at least~$2$, we have $l \geq 2$. By Lemma~\ref{lem:Pc} and \cite[Lemma~4.3]{CH06}, the reflection
subgroup $W(M_i)$ is a direct product of infinite irreducible Coxeter groups. In particular $W(M_i)$ is
center-free by Lemma~\ref{lem:Normalizer}.

Let $M_0 = \Meucl \subset \MF$ be the subset defined after Remark~4.4 in \cite{CH06}. The group $W(M_0)$ is a
direct product of finitely many irreducible affine Coxeter groups of rank $\geq 3$ by
\cite[Proposition~4.7]{CH06}. Moreover, by \cite[Lemma~4.6(ii)]{CH06}, we have either $M_i \cap \Meucl =
\varnothing $ or $M_i \cap \Meucl = M_i$ for each $i \in \{1, \dots, l\}$. Thus, without loss of generality, we
may and shall assume that $M_0 = \Meucl = M_{k+1} \cup M_{k+2} \cup \dots \cup M_l$ for some $k$ (with $k=l$ if
$\Meucl = \varnothing$).

Let now $i \leq k$. By \cite[Lemma~4.6(iii)]{CH06}, for each $j \neq i$, the subgroups $W(M_i)$ and $W(M_j)$
centralize each other. Since moreover, both of them are center-free, we obtain
$$W(\MF) \simeq W(M_0) \times \dots \times W(M_k).$$
The fact that this decomposition passes to the parabolic closure follows easily from Lemma~\ref{lem:DirectProd}
and the fact that finitely generated infinite irreducible Coxeter groups have trivial centre.

Finally, if $M_0 = \Meucl$ is empty, then the fact that $k \geq \dim F$ follows from the proof of
\cite[Theorem~5.2]{CH06}, but only the fact that $k\geq 2$ is relevant to our later purposes.
\end{proof}

An immediate consequence is the following:
\begin{cor}\label{cor:W(MF)}
  Assume that condition (i) of Corollary~\ref{thm:A} holds. Then, for each flat $F$ in $X(W, S)$, the group
  $\Pc(W(\MF))$ is a direct product of irreducible affine Coxeter groups.\qed
\end{cor}

For the sake of future references, we also record the following important fact:
\begin{prop}\label{prop:Flat:Residue}
Let $F$ be a flat in $X$ and let $P$ denote the parabolic closure of $W(\MF)$. Then:
\begin{enumerate}
\item Given any residue $R$ whose stabilizer is $P$ and any wall $m$ which separates some point of $F$ to its
projection to $R$, the wall $m$ separates $F$ from $R$ and  the reflection $r_m$ centralizes $P$.

\item The flat $F$ is contained in some residue whose stabilizer is $P$.
\end{enumerate}
\end{prop}
\begin{proof}
We denote by $\pi_{R}$ the nearest point projection onto $R$ (one may use either the combinatorial projection as defined in \cite[\S3.19]{Ti74}, or
the $\CAT0$ orthogonal projection; both play an equivalent role for our present purposes). Let $x \in F$ be any
point. Since walls and half-spaces in $X$ are closed and convex and since $R$ is $P$-invariant, it follows that
the set $\mathscr{M}(x, \pi_{R}(x))$ of all walls which separate $x$ from $\pi_{R}(x)$ intersects
$\mathscr{M}(R)$ trivially. In other words, the geodesic segment joining $x$ to $\pi_{R}(x)$ does not cross any
wall of $\mathscr{M}(R)$. Therefore, since $\MF \subset \mathscr{M}(R)$, any wall in $\mathscr{M}(x,
\pi_{R}(x))$ separates $F$ from $R$ and, hence, meets every element of $\mathscr{M}(F)$.

Pick an element $m \in \mathscr{M}(x, \pi_{R}(x))$ and let $\mu \in \MF$ be any wall. If the reflections $r_m$
and $r_\mu$ do not commute, then the wall $r_m(\mu)$ is distinct from $m$. Furthermore $r_m(\mu)$ also separates
a point of $F$ to its projection to $R$. By the above, this walls therefore separates $F$ from $R$ and hence, it
belongs to $\mathscr{M}(x, \pi_{R}(x))$. Since the latter set of walls is finite, this shows that the subset of
$\MF$ consisting of all those walls $\mu $ such that $r_m$ does not commute with $r_\mu$ is finite.

By \cite[Lemma~4.3]{CH06}, given any wall $\mu_0 \in \MF$, there exist infinitely many pairs of (pairwise
distinct) walls $\{\mu, \mu' \} \subset \MF$ such that $\mu_0 $ separates $\mu$ from $\mu'$. The preceding
paragraph therefore shows that, given $\mu_0$, we may choose $\mu$ and $\mu'$ in such a way that $r_m$ commutes
with both $r_{\mu}$ and $r_{\mu'}$. Since $r_m$ does not belong to the parabolic subgroup $P$ which, by
definition, contains the parabolic closure of $\la r_\mu, r_{\mu'}\ra$. Therefore Lemma~\ref{lem:debile} implies
that $r_m$ commutes with $r_{\mu_0}$. This implies that $r_m$ centralizes $W(\MF)$ and, hence, normalizes $P$.
Once again, since $r_m$ does not belong to $P$, it follows from Lemma~\ref{lem:Normalizer} that $r_m$
centralizes $P$, thereby establishing (i).

For assertion (ii), choose $R$ amongst the residues whose stabilizer is $P$ in such a way that it minimizes the
distance to $F$. If some point of $F$ does not belong to $R$, there exists a wall which separates that point
from its projection to $R$. By (i) this walls separates $F$ from $R$ and is perpendicular to every wall of $R$.
Transforming $R$ by the reflection through that wall, we obtain another residue whose stabilizer is $P$, but
closer to $F$. This contradicts the minimality assumption made on $R$.
\end{proof}

%\begin{cor}
%Any maximal flat of $X$ can be approximated by periodic flats.
%\end{cor}

\section{Relative  hyperbolicity}

\subsection{Peripheral subgroups are parabolic}%

The purpose of this section is to prove Theorem~\ref{thm:periph}. We will need a subsidiary result on Coxeter
groups. In order to state it properly, we make use of some additional terminology which we now introduce.

Given an element $w \in W$ and a half-space $\mathcal{H}$ of the Cayley graph $\mathrm{Cay}(W, S)$ (or of the
Davis complex $X(W, S)$), we say that $\mathcal{H}$ is \textbf{$w$-essential} is $w.\mathcal{H} \subsetneq
\mathcal{H}$ or $w\inv.\mathcal{H} \subsetneq H$. Notice that an element $w \in W$ admits a $w$-essential
half-space if and only if it has infinite order.

The reflection of $W$ associated to $\mathcal{H}$ is denoted by $r_\mathcal{H}$.

\begin{lem}\label{lem:essential}
Let $H < W$ be a subgroup. Suppose that for any $w \in H$ and any $w$-essential half-space $\mathcal{H}$, the
reflection $r_{\mathcal H}$ belongs to $H$. Then $H$ contains a parabolic subgroup of $W$ as a normal subgroup
of finite index.
\end{lem}
\begin{proof}
Let $P < H$ be the subgroup of $H$ generated by all reflections $r_\mathcal{H}$ associated to a $w$-essential
half-space $\mathcal{H}$ for some element $w \in H$. Thus $P$ is a reflection subgroup of $W$ contained in $H$.
In particular $P$ is itself a Coxeter group, see Lemma~\ref{lem:Pc}.

A crucial point, which follows from \cite[Th.~5.8.1]{Kra94} and \cite[Lem.~5.3]{CH06}, is that $W$ admits a
finite index subgroup $W'$ such that for all $w \in W'$, we have
$$\Pc(w) = \la r_{\mathcal H} \; | \; \mathcal H \text{ is a } w\text{-essential half-space}\ra.$$
In particular $W' \cap H$ is contained in $P$ and hence $P$ has finite index in $H$.

We now choose $w \in P$ in such a way that in the Coxeter group $P$, the parabolic closure $\Pc_P(w)$ of $w$
relative to $P$ is the whole $P$. Such an element $w$ always exists, see \cite[Cor.~3.3]{CapraceFujiwara}. Let
also $P'$ denote the parabolic subgroup of $W$ generated by all those reflections $r_{\mathcal H}$ such that
$\mathcal H$ is a $w$-essential half-space. By the definition of $P$, we have $P' \subset P$. By the property
recalled in the preceding paragraph, the group $P' $ has finite index in $\Pc_P(w) = P$. If follows that $P'$ is
a parabolic subgroup of $W$ which is contained as a finite index subgroup in $P$. In particular $P'$ is a
parabolic subgroup of $W$ which is contained as a finite index subgroup in $H$. Since any intersection of
parabolic subgroups is parabolic, and since $P'$ has finitely many conjugates in $H$, the desired result
follows.
\end{proof}

\begin{proof}[Proof of Theorem~\ref{thm:periph}]
Consider the graph $K$ with vertex set $V = W \cup (\bigcup_{i=1}^m W/H_i)$ and edge set defined as follows. Two
elements of $W$ are joined by an edge if their quotient is an element of $S$; an element $w \in W$ is joined to
a coset $v \in W/H_i$ if and only if $w \in v$. Then $K$ is a connected graph on which $W$ acts by
automorphisms, and containing the Cayley graph of $(W, S)$ as an induced subgraph. By relative hyperbolicity,
this graph is hyperbolic. Furthermore for any $n$ and any two vertices $x, y$, the collection of arcs of length
$n$ joining $x$ to $y$ is finite. In other words, the graph $K$ is \emph{fine} in the terminology of
\cite{Bowditch_relhyp}. Upon adding an edge between any two vertices at distance $2$ in the Cayley graph of
$(W,S)$, we may assume that the graph $K$ has no cut-vertex.

We now apply Bowditch's construction of a hyperbolic $2$-complex $X(K)$ starting from $K$, see Theorem~3.8 (and
also Lemmas~2.5 and~3.3) in \cite{Bowditch_relhyp}. As explained in \emph{loc.~cit.} the action of $W$ on the
boundary of the space $X(K)$ is a geometrically finite convergence action, and the peripheral subgroups (namely
the $W$-conjugates of the $H_i$'s) are the stabilizers of the parabolic points in $\partial X(K)$. In particular
any infinite order element $h \in H_i$ acts as a parabolic element on $X(K)$; it has a unique fixed point $\xi_i
\in
\partial X(K)$ and the limit set of $\la h\ra$ is precisely $\{\xi_i\}$.

\medskip
Let now $r \in W$ be a reflection. Then $r$ acts on the Cayley graph $\mathrm{Cay}(W, S)$ as a reflection.
Clearly $r$ also acts as a reflection on $K$, in the sense that is interchanges two non-empty convex subgraphs
whose union is the whole $K$. It follows from the construction of $X(K)$ that $r$ acts on $X(K)$ as a
\emph{quasi-reflection}: the two half-spaces $\mathcal{H}, \mathcal{H}'$ of $K$ which are interchanged by $r$
correspond in $X(K)$ to two subcomplexes $X(\mathcal{H}), X(\mathcal{H}')$ interchanged by $r$. It follows
immediately that these two subcomplexes are quasi-convex; the fixed point set of $\xi_i$ at infinity thus
coincides with $\partial X(\mathcal{H}) \cap \partial X(\mathcal{H}')$.

Let now $w \in H_i$ be an infinite order element and let $\xi_i \in \partial X(K)$ be the parabolic point fixed
by $H_i$. Let also $\mathcal{H}$ be a half-space of $K$ such that $w(\mathcal{H}) \subsetneq \mathcal{H}$ and
denote by $\mathcal{H}'$ the complementary half-space. Then, for $n>0$ large enough we have $w(X(\mathcal{H}))
\subset X(\mathcal{H})$. Since $w^n.x $ tends to $\xi_i$ with $w \to \infty$ for each $x \in X(K)$, we deduce
that $\xi_i$ belongs to $\partial X(\mathcal{H})$. Applying the same argument to $h\inv$, we deduce on the other
hand that $\xi_i$ belongs to $\partial X(\mathcal{H}')$. Thus $\xi_i \in
\partial X(\mathcal{H}) \cap \partial X(\mathcal{H}')$. It follows that the reflection $r_{\mathcal H}$ fixes $\xi_i$ and,
hence, that $r_{\mathcal H}$ belongs to $H_i$.

By Lemma~\ref{lem:essential}, each peripheral subgroup $H_i$ contains a parabolic subgroup $P_i$ as a finite
index normal subgroup. Since $P_i$ has finite index in $H_i$ it follows that $\xi_i$ is the unique fixed point
of $P_i$ in $\partial X(K)$. In particular the normalizer of $P_i$ in $W$ fixes $\xi_i$. It follows that $H_i =
\norma_W(P_i)$. Upon replacing $P_i$ by a finite index subgroup, we may assume that each irreducible component
of $P_i$ is infinite. It then follows from Lemma~\ref{lem:Normalizer} that $\norma_W(P_i) = H_i$ is itself
parabolic.
\end{proof}

\subsection{Relatively hyperbolic Coxeter groups}

The aim of this section is the proof of Theorem~\ref{thm:RH}. We treat the different implications successively.

\medskip \noindent (ii) $\Rightarrow$ (iii).

Follows from \cite[Theorems~A.0.1 and~A.0.3]{HruskaKleiner_IsolatedFlats}.

\medskip \noindent (ii) $\Rightarrow$ (i).

If (RH1) fails, then $W$ contains a free abelian subgroup which is not contained in any element of
$\mathscr{P}$. If (RH2) fails, then some infinite order element of $W$ is contained in two distinct elements of
$\mathscr{P}$. Therefore, relative hyperbolicity of $W$ with respect to $\mathscr{P}$ implies that (RH1) and
(RH2) both hold.

\medskip \noindent (i) $\Rightarrow$ (iii).

We start with a trivial observation. By condition (RH2), for each irreducible affine subset $J \subset S$, there
is a unique $ J_0 \in \mathscr{T}$ containing $J$. Similarly for each irreducible non-spherical subset $J
\subset S$ such that  $J^\perp$ is non-spherical, there is a unique $ J_0 \in \mathscr{T}$ containing $J$.

For each $J \in \mathscr{T}$, we choose a residue of type $J$ in $X$, which we denote by $R_J$. We define
$\mathscr{F}$ to be the set of all residues of the form $w.R_J$ with $J \in \mathscr{T}$ and $w$ runs over a set
of coset representatives of $\norma_W(W_J)$ in $W$. Note that $\mathscr{T}$, and hence $\mathscr{F}$, is
non-empty unless $W$ is Gromov hyperbolic; of course, we may and shall assume without loss of generality that $W$
is not Gromov hyperbolic.

\smallskip
Let now $F$ be a flat in $X$. Up to replacing it by a conjugate, we may -- and shall -- assume without loss of
generality that $\Pc(W(\MF))$ is standard. Let $I \subset S$ be such that $\Pc(W(\MF)) = W_I$ and $ I_0$ be the
unique element of $\mathscr{T}$ containing $I$. Let $F_0 \in \mathscr{F}$ be the $W_{I_0}$-invariant residue
belonging to $\mathscr{F}$. By Proposition~\ref{prop:Flat:Residue},  any wall $m$ separating a point of $F$ from
its projection to $F_0$ actually separates $F$ from $F_0$; furthermore the reflection $r_m$ through $m$
centralizes $W_{I_0}$. Let $M$ denote the set consisting of all these walls.

Recall that $W(\mathscr{M}(F)) = W_{I}$ is a parabolic subgroup which is a direct product of irreducible
non-spherical subgroups. By (i) and Lemma~\ref{lem:Normalizer}, and since $I_0$ is the unique element of
$\mathscr{T}$ containing $I$, it follows that the centralizer of $W_I$ is contained in the centralizer of
$W_{I_0}$. Furthermore, by the definition of $\mathscr{T}$, the latter centralizer is a finite extension of
$W_{I_0} = W(\mathscr{M}(F_0))$. Since the walls in $M$ may not belong to $\mathscr{M}(F_0)$, it finally follows
that $W(M)$ centralizes $W(\mathscr{M}(F_0))$. Hence $W(M)$ is finite, and so is $\Pc(W(M))$ by
Lemma~\ref{lem:Pc}. In particular, the cardinality of $M$ is bounded above by the maximal number of reflections
in a finite standard parabolic subgroup. This shows that the combinatorial distance from any point $x \in F$ to
$F_0$ is uniformly bounded. Therefore, condition (A) is satisfied.

Now we prove (B). Let $F, F'  \in \mathscr{F}$ be residues such that $\mathcal{N}_r(F) \cap \mathcal{N}_r(F')$
is unbounded for some $r >0$. Then the visual boundaries $\partial_\infty(F)$ and $\partial_\infty(F')$ have a
common point. In other words, there exists a geodesic ray $\rho \subset F$ and $\rho' \subset F'$ such that
$\rho$ and $\rho'$ are at bounded Hausdorff distance. Let $\mathscr{M}(\rho)$ be the set of walls which separate
two points of $\rho$. Since $\mathscr{M}(\rho)$ is infinite whereas for any $x \in \rho$, the set of walls which
separate $x$ from $\rho'$ is uniformly bounded, it follows that $\mathscr{M}(\rho) \cap \mathscr{M}(\rho')$ is
infinite. Therefore $\mathscr{M}(\rho) \cap \mathscr{M}(\rho')$ contains two walls $m, m'$ which do not meet
\cite[Lemma~13]{Ca06}.

Denote by $P, P'$ the respective stabilizers of $F, F'$ in $W$. Notice that $P$ and $P'$ are parabolic subgroups
whose reflections consist of the sets $\mathscr{M}(F)$ and $\mathscr{M}(F')$ respectively. The preceding
paragraph shows that $P$ and $P'$ share  a common infinite dihedral subgroup. By (RH2), this implies that $P$
and $P'$ coincide. In view of the definition of $\mathscr{F}$, we deduce that $F$ and $F'$ must coincide.

This shows that for any two distinct $F, F' \in \mathscr{F}$ and each $r > 0$, the set $\mathcal{N}_r(F) \cap
\mathcal{N}_r(F')$ is bounded. The fact that its diameter depends only on $r$, but not on the specific choice of
$F$ and $F'$, follows from the cocompactness of the $W$--action on $X$. Hence (B) holds.

\medskip \noindent (iv) $\Rightarrow$ (iii).

Clear since the Davis complex $X(W, S)$ is a (thin) building of type $(W, S)$.

\medskip \noindent (iii) $\Rightarrow$ (iv).

Let $Y$ be a building of type $(W,S)$ and $\mathscr{F}$ be the set of all $J$-residues of $Y$ with $J \in
\mathscr{T}$. Furthermore, given an apartment $A$ in $Y$, we set
$$\mathscr{F}_A = \{A \cap F \; | \; F \in \mathscr{F}, A \cap F \neq \varnothing\}.$$
Since (iii) holds and since an apartment in $Y$ is nothing but an isometrically embedded copy of the Davis
complex $X(W, S)$, it follows that the elements of $\mathscr{F}_A$ are isolated in $A$. Moreover, the constant
$D$ which appears in condition (A) depends only on $(W, S)$.

Let $F$ be a flat in $Y$. Then $F$ is contained in an apartment $A$ by \cite[Theorem~6.3]{CH06}. Therefore
condition (A) holds for $Y$ since the elements of $\mathscr{F}_A$ are isolated in $A$.

Let now $J, J' \in \mathscr{T}$ and $F, F' \in \mathscr{F}$ be residues of type $J$ and $J'$ respectively.
Assume that $\mathcal{N}_r(F) \cap \mathcal{N}_r(F')$ is unbounded for some $r > 0$. Let $A$ be an apartment
contained a chamber $c$ of $F$ and let $c'$ be the combinatorial projection of $c$ onto $F'$ (see
\cite[\S3.19]{Ti74} ). We denote by $\rho_{c, A}$ the combinatorial retraction of $Y$ onto $A$ centered at $c$.
Recall that this maps any chamber $x$ of $Y$ to the unique chamber $x'$ of $A$ such that $\delta_Y(c,
x')=\delta_Y(c, x)$, where $\delta_Y : \Ch(Y) \times \Ch(Y) \to W$ denotes the Weyl distance.

By assumption, there exists an unbounded sequence $c' =c'_0, c'_1, \dots$ of chambers of $F'$ such that $c'_n$
lies at uniformly bounded distance from $F$. Since combinatorial retractions do not increase distances and since
$\rho_{c, A}$ maps any chamber in $F$ to a chamber in $A \cap F$, it follows that the sequence $(x'_n)$ defined
by $x'_n = \rho_{c, A}(c'_n)$ lies at uniformly bounded distance from $A \cap F$. Furthermore, by a standard
property of the combinatorial projection, namely the \emph{gate property} (see \cite[Ch.~3]{Ti74}), for each chamber $x' \in F'$ there exists a minimal gallery joining
$c$ to $x'$ via $c'$. Therefore, it follows that for each $n$, the chamber $x'_n$ lies in the $J'$-residue
containing $x'_0$, say $F''$. This shows in particular that $\mathcal{N}_r(A \cap F) \cap \mathcal{N}_r(A \cap
F'')$ is unbounded. Since $F''$ is a residue of type $J'$ and since the residues whose type belong to
$\mathscr{T}$ are isolated in $A$ by assumption, it follows that $A \cap F = A \cap F''$ and hence $F= F''$ and
$J=J'$. In particular we obtain $x'_0 \in F$ because $c \in  F$ and $x'_0 \in F''$. Since $\delta_Y(c, x'_0) =
\delta_Y(c, c')$, we deduce that $c' \in F$, whence $F = F'$ since $J=J'$.

This shows that for any two distinct $F, F' \in \mathscr{F}$ and each $r > 0$, the set $\mathcal{N}_r(F) \cap
\mathcal{N}_r(F')$ is bounded. The fact that its diameter depends only on $r$, but not on the specific choice of
$F$ and $F'$, follows from the corresponding fact for apartments in $Y$. Hence (B) holds.\qed

\begin{proof}[Proof of Corollary~\ref{cor:MaxParab}]
It is immediate to check that the set $\mathscr{T}$ satisfies (RH1) and (RH2).
\end{proof}

\subsection{Isolated flats}

\begin{lem}\label{lem:RH:aff}
The following conditions are equivalent:
\begin{enumerate}
\item The collection $\mathscr{T}$ of maximal Euclidean subsets of $S$ satisfies (RH1) and (RH2).

\item For all non-spherical $J_1, J_2 \subset S$ such that $[J_1, J_2]=1$, the group $\la J_1 \cup J_2 \ra$
    is virtually abelian.

\item For each minimal hyperbolic $J \subset S$, the set $J^\perp$ is spherical.
\end{enumerate}
\end{lem}
\begin{proof} The main point is that, given a Coxeter system $(W, S)$, it is well known the group $W$ is
virtually abelian if and only if it is a direct product of finite and affine Coxeter groups, i.e. if $S$ is
Euclidean (see e.g. \cite{MarVin00}).

\noindent (i) $\Rightarrow$ (iii).

Let $\mathscr{T}$ be a collection of subsets of $S$ satisfying (RH1) and (RH2). If (iii) fails, then there exist
a minimal hyperbolic subset $J$ and a non-spherical irreducible subset $I \subset J^\perp$. By (RH1) there
exists $K \in \mathscr{T}$ such that $I \cup J \subset K$. Then $\la K \ra$ is not virtually abelian since it
contains $\la J\ra$, hence (i) fails as well.

\noindent (iii) $\Leftrightarrow$ (ii).

If (ii) fails then there exists a non-spherical and non-affine subset $J\subset S$ such that $J^\perp$ is
non-spherical. Now any minimal non-spherical and non-affine subset $I$ of $J$ is minimal hyperbolic, and since
$I \subset J$ we have $I^\perp \supset J^\perp$. Thus (iii) fails as well.

\noindent (ii) $\Rightarrow$ (i).

The condition (ii) clearly implies that for each irreducible non-spherical subset $J$, either $J$ is affine and
$J \cup J^\perp$ is Euclidean or $J$ is non-affine and $J^\perp$ is spherical. In particular, every irreducible
affine subset is contained in a \emph{unique} maximal Euclidean subset. In other words the collection
$\mathscr{T}$ of maximal Euclidean subsets of $S$ satisfies (RH2). Moreover (RH1) clearly holds as well.
\end{proof}

\begin{proof}[Proof of Corollary~\ref{thm:A}]
In view of Theorem~\ref{thm:RH}, Lemma~\ref{lem:RH:aff} and \cite[Th.~1.2.1]{HruskaKleiner_IsolatedFlats}, it is
enough to prove that $W$ is relatively hyperbolic with respect to its maximal virtually abelian subgroups if and
only if it is relatively hyperbolic with respect to its maximal parabolic subgroups of Euclidean type. Since any
parabolic subgroup of Euclidean type is virtually abelian, the `if' part is clear. Conversely, assume that $W$
is relatively hyperbolic with respect to its maximal virtually abelian subgroups. Then conditions (i) and (ii)
hold. In view of \cite[Prop.~3.2]{CM05}, this implies that the parabolic closure of any virtually abelian
subgroup of rank $\geq 2$ is of Euclidean type. In particular, if $A < W$ is a maximal virtually abelian
subgroup, then $A = \Pc(A)$.
\end{proof}

\begin{proof}[Proof of Corollary~\ref{cor:B}]
Assertions (i) and (iii) follow from \cite[Theorem~A.0.1]{HruskaKleiner_IsolatedFlats}. For (iv), one applies
\cite[Proposition~5.4]{DrutuSapir}; one needs the fact that any asymptotic cone of a Euclidean building is
itself a Euclidean building: this is established in \cite[Theorem~1.2.1]{KleinerLeeb}. The fact that $W$ is
biautomatic follows from \cite[Theorem~1.2.2(5)]{HruskaKleiner_IsolatedFlats}. The biautomaticity of $\Gamma$
can then be deduced either directly from \cite{Swiatk} or from (i) together with \cite{Rebbechi} and the fact
that cocompact lattices of Euclidean buildings are biautomatic by \cite{Swiatk}.
\end{proof}

\newpage
\appendix

\pagestyle{empty}

\begin{center}
\textbf{\Large{Erratum to `Buildings with isolated subspaces \\ and relatively hyperbolic Coxeter groups'}}

\medskip
 Pierre-Emmanuel Caprace

\medskip
December 2013
\end{center}

\bigskip

%\author{Pierre-Emmanuel Caprace, Colin D. Reid and George A. Willis}

%\author[1]{Pierre-Emmanuel Caprace}
%\ead{pe.caprace@uclouvain.be}

%\affil[1]{Universit\'e catholique de Louvain, IRMP, Chemin du Cyclotron 2, bte L7.01.02, 1348 Louvain-la-Neuve, Belgique}
%\date{November 2013}

%\maketitle

The goal of this note is to correct two independent errors in \cite{Cap_RH}, respectively in Theorems~A and~B from loc. cit. I am indebted to Alessandro Sisto, who pointed them out to me.  Those corrections affect neither the characterization of toral relatively hyperbolic Coxeter groups (Cor.~D and~E from  \cite{Cap_RH}), nor the other intermediate results from the original paper. 

We keep the notation and terminology from loc. cit. Moreover all Coxeter groups under consideration are assumed to be finitely generated. 
The first correction concerns Theorem~A: its assertions (ii), (iii), (iv)  are indeed equivalent, but a third condition (RH3) has to be added to (RH1) and (RH2) in assertion (i), as in the following reformulation. 

\renewcommand{\thethmintro}{\Alph{thmintro}$'$}
\setcounter{thmintro}{0}
\begin{thmintro}\label{thm:A}
Let $(W, S)$ be a Coxeter system and $\mathscr T$ be a set of subsets of the Coxeter generating set $S$. Then $W$ is hyperbolic relative to $\mathscr P = \{W_J \; | \; J \in \mathscr T\}$ if and only if the following three conditions hold:
\begin{description}
\item [\textbf{(RH1)}] For each irreducible affine subset $J \subset S$ of cardinality~$\geq 3$, there exists $K \in \mathscr T$ such that $J \subset K$. Similarly, for each pair of irreducible non-spherical subsets $J_1, J_2 \subset S$ with $[J_1, J_2]=1$, there exists $K \in \mathscr T$ such that $J_1 \cup J_2 \subset K$. 

\item [\textbf{(RH2)}] For all $K_1, K_2 \in \mathscr T$ with $K_1 \neq K_2$, the intersection $K_1 \cap K_2$ is spherical. 

\item [\textbf{(RH3)}] For each $K \in \mathscr T$ and each irreducible non-spherical $J \subset K$, we have $J^\perp \subset K$.
\end{description}
\end{thmintro}

\begin{proof}
The necessity of (RH1) and (RH2) is established in \cite{Cap_RH}. The condition (RH3) is also necessary, as pointed out by Alessandro Sisto: if there is a reflection $s \in S$ and a set $K \in \mathscr T$ such that $s \not \in K $ and $s$ commutes with an irreducible non-spherical subset $J \subset K$, then the cosets $W_K$ and $sW_K$ of the parabolic subgroup $W_K$ are distinct, but the intersection of their respective $1$-neighbourhoods in the Cayley graph of $(W, S)$ is unbounded, since it contains $W_J$. This contradicts that $W$ is hyperbolic relative to $\mathscr P$. 

Assume conversely that (RH1), (RH2) and (RH3) hold. As in \cite{Cap_RH}, we need to show that the set $\mathscr F$, consisting of all residues of the Davis complex of $(W, S)$ whose type belongs to $\mathscr T$, satisfies the isolation conditions (A) and (B) from loc. cit. The arguments given  there show that (RH1) is sufficient to ensure that (A) holds. Moreover it is shown that if $\mathscr F$ does not satisfy (B), then there exists two distinct residues $F, F' \in \mathscr F$ whose respective stabilisers   $P, P'$, which are parabolic subgroups of $W$, share a common infinite dihedral reflection subgroup. The mistake in \cite{Cap_RH} lies in the sentence: `By (RH2), this implies that $P$ and $P'$ coincide.' The corrected argument, which requires also invoking (RH3), goes as follows. We may write $P = g W_{K} g\inv$ and $P' = g' W_{K'} (g')\inv$ for some $K, K' \in \mathscr T$ and $g, g' \in W$. Since $P \cap P'$ contains an infinite dihedral reflection subgroup, it also contains the parabolic closure of that subgroup, say $Q$, which is of irreducible non-spherical type by \cite[Lem.~2.1]{Cap_RH}. Therefore there is an  irreducible non-spherical subset $J \subset K$ (resp. $J' \subset K'$) such that $Q $ is conjugate to $gW_Jg\inv$ in $P$ (resp. to $g'W_{J'}(g')\inv$ in $P'$). It follows that $W_J$ is conjugate to $W_{J'}$ and, hence, that   $J$ and $J'$ are conjugate in $W$. By \cite[Prop.~5.5]{Deodhar}, it follows that $J = J'$, so that   $K = K'$ by (RH2). In particular $P$ and $P'$ are conjugate. Let $p \in P$ be an element which conjugates $g W_J g\inv$ to $Q$.   Upon replacing $g$ by    $pg$, we may assume that $Q= gW_J g\inv$. Similarly we may assume that $Q = g' W_J (g')\inv$. It follows that $g\inv g'$ normalises $W_J$. By \cite[Prop.~5.5]{Deodhar}, the normaliser of $W_J$ coincides with $W_{J \cup J^\perp}$, and is thus contained in $W_K$ by (RH3). Hence $g\inv g'$ normalizes $W_K$, so that $P = P'$. Condition (RH3) together with \cite[Prop.~2.1]{BH} and \cite[Prop.~5.5]{Deodhar} also implies that $P$ is self-normalising, which implies that there is a unique residue in the Davis complex, whose full stabiliser is $P$. We deduce that $F= F'$, a contradiction. This confirms that (B) holds. % and   the proof of Theorem~\ref{thm:A} is complete.  
\end{proof}

We next remark that Corollaries~D and E from \cite{Cap_RH} are not affected by the above correction: indeed, in the respective settings of those corollaries, the condition (RH3) holds automatically. In Corollary~C, for all three conditions (RH1)--(RH3) to be satisfied, the definition of $\mathscr T$ has to be adapted as follows:
$$\mathscr T = \big\{ S \setminus \{s_0\} \big\} \cup \big\{ J \cup J^\perp \; | \; J \text{ is irreducible affine of cardinality $\geq 3$ and contains } s_0\big\}.
$$

\medskip
We now turn to the second error, which lies in Theorem~B from \cite{Cap_RH}. The purpose of that statement was to answer the following question: \emph{assuming that $W$ is hyperbolic with respect to some peripheral subgroups $H_1, \dots, H_m$, can one relate those peripheral  subgroups to the parabolic subgroups of $W$ (in the usual Coxeter group theoretic sense)?} Theorem~B asserted that those peripheral subgroups are always parabolic in the Coxeter group theoretic sense. This is   not true in general: indeed,  any Gromov hyperbolic group is also relatively hyperbolic with respect to any malnormal collection of quasi-convex subgroups, see \cite[Th.~7.11]{Bowditch}. Therefore, even if $W$ is Gromov hyperbolic, one can always make it relatively hyperbolic by adding maximal self-normalising cyclic subgroups as peripheral subgroups, and those are not parabolic in the Coxeter sense. The correct statement can be phrased as follows: \emph{if $W$ is relatively hyperbolic with respect to some peripheral subgroups  $H_1, \dots, H_m$, then it is also relatively hyperbolic with respect to a (possibly empty) collection of Coxeter-parabolic subgroups $P_1, \dots, P_k$, and moreover, each $P_i$ is conjugate to a subgroup of some $H_j$}. In particular every Coxeter group admits a canonical, minimal, relatively hyperbolic structure, whose peripheral subgroups are indeed parabolic in the Coxeter group theoretic sense. The latter result has been obtained in a joint work with Jason Behrstock, Mark Hagen and Alessandro Sisto. In that work, we also provide various characterizations of the canonical parabolic subgroups $P_1, \dots, P_k$, and describe necessary and sufficient conditions on a Coxeter presentation of $W$ ensuring that $W$ is not   relatively hyperbolic with respect to any collection of proper subgroups. Those results appear  in the Appendix to \cite{BHS}.

\end{document}